\documentclass[11pt]{amsart}
\usepackage{amscd,amssymb}
\usepackage[dvips]{graphics}
\usepackage{tabularx}
\input diagrams
\diagramstyle[PostScript=dvips]

\newcommand{\scheme}[1]{\mathcal {#1}}

\newcommand{\mbar}{\M}

\newcommand{\mgnbar}{\overline{\scheme{M}}_{g,n}}
\newcommand{\mgnrmbar} {\mgnbar^{1/r,\bm}}
\newcommand{\mgnrbar} {\mgnbar^{1/r}}

\newcommand{\cgn}{\scheme{C}_{g,n}}
\newcommand{\cgnr}{\cgn^{1/r}}

\newcommand{\pic}{\text{Pic\,}}        

\newcommand{\tensor}{\otimes}

\newcommand{\rk}{\operatorname{rk}}

\newcommand{\ck}{{\mathcal K}}

\newcommand{\ce}{{\mathcal E}}

\newcommand{\co}{{\mathcal O}}

\newcommand{\ch}{{\mathcal H}}

\newcommand{\sheafhom}{\ch\kern-.15em om}

\newcommand{\ba}{\mathbf{a}}
\newcommand{\bdelta}{\boldsymbol{\delta}}

\newcommand{\bm}{\mathbf{m}}          %
\newcommand{\bmt}{\widetilde{\mathbf{m}}}          %

\newcommand{\bt}{\mathbf{t}}
\newcommand{\bTt}{\tilde{\bt}}
\newcommand{\btau}{\boldsymbol{\tau}}
\newcommand{\btaut}{\widetilde{\btau}}

\newcommand{\bx}{\mathbf{x}}          
\newcommand{\cft}{CohFT}              
\newcommand{\chit}{\tilde{\chi}}      %
\newcommand{\cor}[1]{\langle\,{#1}\,\rangle}  
\newcommand{\cort}[1]{\langle\,{#1}\,\rangle}  
\newcommand{\cv}{c^{1/r}}            

\newcommand{\KdV}{\mathrm{KdV}}
\newcommand{\mt}{{\tilde{m}}}         
\newcommand{\M}{\overline{\MM}}       
\newcommand{\MM}{\scheme{M}}          

\newcommand{\nq}{{\mathbb{Q}}}        
\newcommand{\psit}{\tilde{\psi}}     %

\newcommand{\taut}{\tilde{\tau}}
\newcommand{\Tt}{\tilde{t}}

\def\lift#1#2{
  \dimen0 = \unitlength
  \multiply\dimen0 by #1 \divide \dimen0 by 2
  \dimen1 = \dimen0
  \multiply \dimen1 by 7 \divide \dimen1 by 10
  \raise\dimen1
     \hbox{\hskip 0.3cm ${\vbox to \dimen0{}}$ \enspace #2}}

\newtheorem{thm}{Theorem}[section]
\newtheorem{introdax}{Descent Axiom \ref{ax:dax}}%
\newtheorem{dax}[thm]{Descent Axiom}
\newtheorem{vax}[thm]{Vanishing Axiom}
\newtheorem{lm}[thm]{Lemma}
\newtheorem{prop}[thm]{Proposition}
\newtheorem{crl}[thm]{Corollary}

\theoremstyle{definition}

\newtheorem{rem}[thm]{Remark}

\newtheorem{df}[thm]{Definition}

\theoremstyle{remark}
 
\newtheorem{ack}{Acknowledgments}

\begin{document}
\addtocounter{section}{-1}

\title[Gravitational Descendants]
{Gravitational Descendants and the Moduli Space of Higher Spin
Curves} 
\subjclass{Primary: 14N35, 53D45. Secondary: 14H10} 

\author
[T. J. Jarvis]{Tyler J. Jarvis}
\address
{Department of Mathematics, Brigham Young University, Provo, UT
84602, USA} \email{jarvis@math.byu.edu}
\thanks{Research of the first author was partially supported by NSA
grant MDA904-99-1-0039}

\author
[T. Kimura]{Takashi Kimura}
\address
{Department of Mathematics, 111 Cummington Street, Boston
University, Boston, MA 02215, USA} \email{kimura@math.bu.edu}
\thanks{Research of the second author was partially supported by NSF grant
  number DMS-9803427}

\author
[A. Vaintrob]{Arkady Vaintrob\vskip 0.2cm }
\address
{Department of Mathematics, University of Oregon, Eugene, OR
97403, USA} \email{vaintrob@math.uoregon.edu}

\date{}

\begin{abstract}

The purpose of this note is introduce a new axiom (called the
Descent Axiom) in the theory of $r$-spin cohomological field
theories. This axiom explains the origin of gravitational
descendants in this theory. Furthermore, the Descent Axiom
immediately implies the Vanishing Axiom, explicating the latter
(which has no \emph{a priori} analog in the theory of 
Gromov-Witten invariants), in terms of the multiplicativity of the 
virtual class.  We prove that the Descent Axiom holds in the 
convex case, and consequently in genus zero. 
\end{abstract}

\maketitle

\section{Introduction}
\label{intro}

In \cite{JKV}, the notion of an $r$-spin cohomological field theory was
introduced for each integer $r\geq 2$. This is a particular realization of a
cohomological field theory (\cft) in the sense of Kontsevich-Manin
\cite{KM1}.  Its construction was motivated by drawing an analogy \cite{JKV2}
with the Gromov-Witten invariants of a smooth, projective variety $V$.

The analog of the moduli space $\M_{g,n}(V)$ of stable maps into
$V$ is $\M_{g,n}^{1/r}$, the moduli space of stable $r$-spin
curves \cite{J,J2}. The moduli space $\M_{g,n}^{1/r}$ is the
disjoint union of $\M_{g,n}^{1/r,\bm}$ where $\bm =
(m_1,\ldots,m_n)$ and $m_i = 0,\ldots,r-1$ for all $i=1,\ldots,n$.
Suppose that $r\geq 2$ is prime. For any $n$-tuple of integers
$\bm$, $\M_{g,n}^{1/r,\bm}$ is a compactification of the moduli
space of (connected) Riemann surfaces $\Sigma$ of genus $g$ with
$n$ marked points,  $(p_1,\ldots,p_n)$, together with a
holomorphic line bundle $\ce$ over the curve whose $r$-th tensor
power is $\ck$, the canonical bundle of $\Sigma$, twisted by
$\co(-\sum_{i=1}^n m_i p_i)$, together with a bundle isomorphism
between $\ce^{\otimes r}$ and $\ck\otimes\co(-\sum_{i=1}^n m_i
p_i)$. When $r$ is not prime then one must include, in addition,
$d$-th roots for all $d$ dividing $r$ and suitable compatibility
morphisms. Furthermore, when the $m_i$ are permitted to range over  
all non-negative integers, any two spaces $\M_{g,n}^{1/r,\bm}$ and 
$\M_{g,n}^{1/r,\bm'}$ are canonically isomorphic  if  the 
$n$-tuples $\bm$ and $\bm'$ are componentwise equal mod $r$. 
Therefore, one may restrict without loss of generality to the 
space $\M_{g,n}^{1/r}$. In particular, the space $\M_{g,n}^{1/r}$ 
is a smooth stack which is a ramified cover of the moduli space of 
stable curves $\M_{g,n}$.

There is an analog  $\cv$ of the cap product of the virtual
fundamental class of $\M_{g,n}(V)$ with the Gromov-Witten classes
on $\M_{g,n}(V)$ in the $r$-spin \cft\, called a \textsl{virtual
class} on $\M_{g,n}^{1/r}$. Following the ideas of Witten
\cite{W,W2}, we proved in \cite{JKV} that if $\cv$ satisfied
certain axioms then it would give rise to a \cft. The class $\cv$
was constructed for all genus when $r=2$ and in genus zero for all
$r\geq 2$.

One can form the large phase space potential $\Phi(\bt)$
associated to the $r$-spin \cft\ which is a generating function
for the correlators (including the gravitational descendants) of
the theory. It was proved in \cite{JKV} that for all $r\geq 2$,
the part of the potential corresponding to genus zero solves (a
semi-classical limit of) the $r$-th Gelfand-Dickey (or $\KdV_r$)
integrable hierarchy following the ideas of Witten \cite{W,W2}.
Witten conjectured that this correspondence should extend to all
genera. When $r=2$, this conjecture reduces to Kontsevich's result
\cite{Ko} relating the usual $\KdV_2$ hierarchy to intersection
numbers on the moduli space $\M_{g,n}$ of stable curves \cite{W3}.

While most of the axioms satisfied by $\cv$ have analogs in the
theory of Gromov-Witten invariants, there is one axiom satisfied
by $\cv$, called the Vanishing Axiom,  which appears to have no
analog in the theory of Gromov-Witten invariants. The Vanishing
Axiom states that $\cv$ vanishes on $\M_{g,n}^{1/r,\bm}$ if $\bm =
(m_1,\ldots,m_n)$ and $m_i = r-1$ for some $i$ and $0\leq m_j\leq
r-1$ for all $j$. Since this axiom does not immediately appear to
have the form of a factorization identity, it does not appear to
have an analog in the theory of Gromov-Witten invariants.

There is another puzzling feature in the $r$-spin \cft. Although 
the class $\cv(\bm)$ is defined on $\M_{g,n}^{1/r,\bm}$ where $\bm 
= (m_1,\ldots,m_n)$ and $m_i = 0,\ldots,r-1$ for all 
$i=1,\ldots,n$, the axioms of the virtual class are such that one 
may extend its definition to $\cv(\bmt)$ in 
$H^\bullet(\M_{g,n}^{1/r})$ for arbitrary  $n$-tuples of 
nonnegative integers $\bmt = (\mt_1,\ldots,\mt_n)$ under the 
identification of $\M_{g,n}^{1/r,\bm}$ with $\M_{g,n}^{1/r,\bm'}$ 
whenever the $n$-tuples of integers $\bm$ and $\bm'$ are 
equivalent (componentwise) mod $r$. Furthermore, it is simple to 
see that the usual construction of the virtual class in genus 
zero, extends straightforwardly to yield a class $\cv(\bmt)$ for  
arbitrary nonnegative $n$-tuples $\bmt$. It is natural to ask what 
role these additional (infinite number of) classes $\cv(\bmt)$ 
play in the $r$-spin \cft.

Finally, there is a broader issue that permeates both the theory
of Gromov-Witten invariants and the $r$-spin \cft. In the case of
$r$-spin \cft\ (the Gromov-Witten case is analogous),  the
gravitational descendants are correlators which are integrals of
products  of tautological $\psi_i$ classes on $\M_{g,n}^{1/r}$
(where $i=1,\ldots,n$) with $\cv$. From a purely geometric
perspective, the appearance of the $\psi_i$ classes in the
correlators is rather ad hoc.  The real question is: \textsl{Where
do these $\psi$ classes come from?}

All of these issues can be simultaneously addressed if one assumes
that $\cv$ satisfies a new axiom called \textsl{the Descent
Axiom}.

\begin{introdax}
Let $r\geq 2$ be an integer, and $\bm = (m_1,\dots,m_n)$ be an
$n$-tuple of integers such that either $m_i \geq 0$ for all
$i=1,\ldots,n$, or there exists an integer  $1\leq j\leq n$ such
that $m_i \geq 0$ for all $i\not=j$ and $m_j = -1$.
Let $\cv(\bm)$ denote the virtual class on $\M_{g,n}^{1/r,\bm}$.
Let $r \bdelta_i$ denote an $n$-tuple of integers which is $r$ in
the $i$th position and zero in all others then
\[
\cv(\bm+ r \bdelta_i) = -\psit_i(\bm) \cv(\bm)
\]
\end{introdax}
\renewcommand{\thedax}{\arabic{dax}}

The Vanishing Axiom can be shown to follow from properties of the
class $\psit_i(\bm_i)$.

The Descent Axiom implies that if one considers a generating
function of intersection numbers obtained by integrating ONLY
$\cv(\bm)$ for \emph{all} $\bm$ with $m_i \geq 0$ against the
fundamental class of $\mgnrbar$ then one obtains, after a simple 
variable redefinition, the usual large phase space potential of 
the $r$-spin \cft! In other words, the Descent Axiom explains 
geometrically why the $\psi$ classes appear in the usual 
definition of gravitational descendants.

We will prove that the Descent Axiom holds in the convex case
(and, hence, in genus zero) using the multiplicativity property of
the top Chern class. The Vanishing Axiom can then be interpreted
as following from a kind of factorization property of $\cv$.

Finally, it is worth pointing out that all of the ideas in this
note generalize to the moduli space of case of Gromov-Witten
invariants through the introduction of the moduli space of stable
$r$-spin maps into $V$ \cite{JKV3}.

In the first section of this note, we briefly review the
properties of the moduli space of stable $r$-spin curves
$\M_{g,n}^{1/r,\bm}$ and the tautological cohomology classes
$\psi_i$ and $\psit_i(\bm)$ associated to tautological line 
bundles over $\M_{g,n}^{1/r,\bm}$. We then recall the Vanishing 
Axiom and prove that the Descent Axiom implies the Vanishing Axiom 
by using a universal relation between these classes.

In the second section, we show that the Descent Axiom implies that
the usual large phase space potential function of an $r$-spin
\cft\ $\Phi(\bt)$, containing both $\cv$ and the $\psi_i$ classes,
agrees (after a simple variable redefinition) with the generating
function for correlators associated to the virtual classes
$\cv(\bm)$ where $\bm=(m_1,\ldots,m_n)$ are arbitrary, nonnegative
$n$-tuples.

In the third section, we recall the construction of the virtual
class $\cv$ in the convex case and prove that the  Descent Axiom
holds in this case.

\begin{ack}
Parts of this note were written while T.K.\ was visiting the
Universit\'{e} de Bourgogne  and A.V.\ was visiting the Institut
des Hautes \'Etudes Scientifiques. We would like to thank these
institutions for their hospitality and support.
\end{ack}

\section{The Moduli Stack of $r$-spin Curves}

We will use the notation and conventions of \cite{JKV}, which we
briefly review here.  Complete details and proofs may be found in
\cite{JKV} and in \cite{J}.

\subsection{Definitions and Basic Properties}
\begin{df}\label{def:prestable}
A \emph{prestable curve} is a reduced, complete, algebraic curve
with at worst nodes as singularities.
\end{df}

\begin{df}\label{def:root}
Let $(X, p_1, \dots, p_n)$ be a prestable, $n$-pointed, algebraic
curve, let $\ck$ be a rank-one, torsion-free sheaf, and let $\bm =
(m_1,\ldots,m_n)$  be a collection of integers. A $d$-th root of
$\ck$ of type $\bm$ on $X$ is a pair $(\ce,c)$ of a rank-one,
torsion-free sheaf $\ce$ and an $\co_X$-module homomorphism $c:
\ce^{\tensor d} \to \ck \tensor \co_X ( - \sum_{i=1}^{n} m_i p_i)$
with the following properties:
\begin{itemize}
\item $d \cdot \deg \ce = \deg \ck - \sum m_i$
\item $c$ is an isomorphism on the locus of $X$ where $\ce$ is
locally free
\item for every point $p \in X$ where $\ce$ is not free, the
length of the cokernel of $b$ at $p$ is $d-1$.
\end{itemize}
\end{df}

Although the condition on the cokernel seems strange, it turns out
to be a very natural one (see \cite{J}). Indeed,
the construction of $r$-spin curves can
also be done \cite{AJ} in terms of line bundles on the ``twisted
curves'' of Abramovich and Vistoli, and in their formulation, the
cokernel condition amounts exactly to the condition that the local
(orbifold) index of the curve at each node or marked point must
divide $d$---a redundant condition in that formulation.

For any $d$-th root $(\ce,c)$ of $\ck$ of type $\mathbf{m}$, and
for any $\mathbf{m}'$ congruent to $\mathbf{m} \pmod d$, we can
construct a unique $d$-th root $(\ce',b')$ of type $\mathbf{m}'$
simply by taking $\ce'=\ce \otimes \co(1/d \sum (m_i-m_i')p_i)$.
Consequently, the moduli of curves with $d$-th roots of a bundle
$\ck$ of type $\mathbf{m}$ is canonically isomorphic to the moduli
of curves with $d$-th roots of type $\mathbf{m}'$.

\begin{df}\label{def:spin-structure}
Let $(X, p_1, \ldots, p_n)$ be a prestable, $n$-pointed curve. An
\emph{$r$-spin structure on $X$ of type $\bm =(m_1, \ldots, m_n)$}
is a pair $(\{\ce_d\}, \{c_{d, d'}\})$ of a set of sheaves and a
set of homomorphisms as follows.  The set of sheaves consists of a
rank-one, torsion-free sheaf $\ce_d$ on $X$ for every divisor $d$
of $r$; and the set of homomorphisms consists of an $\co_X$-module
homomorphism $ c_{d,d'} : \ce^{\tensor d/d'}_{d} \rTo \ce_{d'} $
for every pair of divisors $d',d$ of $r$,  such that $d'$ divides
$d$. These sheaves and homomorphisms must satisfy the following
conditions:
\begin{itemize}
\item The sheaf $\ce_1$ is isomorphic to the canonical (dualizing)
sheaf $\omega_X$.
\item For each divisor $d$ of $r$ and each divisor $d'$ of $d$,
the  homomorphism $c_{d,d'}$  makes $(\ce_d, c_{d,d'})$ into a
$d/d'$-th root of $\ce_{d'}$ of type $\mathbf{m}'$, where
$\mathbf{m}'=(m'_1, \ldots, m_n')$ is the reduction of $\bm$
modulo $d/d'$ (i.e. $0\le m_i' < d/d'$ and $m_i \equiv m_i' \pmod
d/d'$).
\item The homomorphisms $\{c_{d,d'}\}$ are compatible.
That is, for any integers $d$ dividing $e$ dividing $f$ dividing
$r$, we have  $c_{f,d} = c_{e,d} \circ c_{f,e}$.
\end{itemize}
\end{df}

If $r$ is prime, then an $r$-spin structure is simply an $r$-th
root of $\omega_X$.  Even when $d$ is not prime, if the root
$\ce_d$ is locally free, then for every divisor $d'$ of $d$, the
sheaf $\ce_{d'}$ is uniquely determined, up to an automorphism of
$\ce_{d'}$. In particular, if $\bm'$ satisfies the conditions
$\bm' \equiv \bm \pmod {d'}$ and $0 \leq m'_i <d'$,  then the
sheaf $\ce_{d'}$ is isomorphic to $\ce_{d}^{\tensor d/d'}\otimes
         \co\left(\frac{1}{d'} \sum (m_i -m'_i)p_i\right)$.

Every $r$-spin structure on a smooth curve $X$ is determined, up
to isomorphism, by a choice of a line bundle $\ce_r$, such that
$\ce^{\tensor r}_r \cong \omega_X (-\sum m_i p_i)$.  In
particular, if $X$ has no automorphisms, then the set of
isomorphism classes of $r$-spin structures (if non-empty) of type
$\mathbf{m}$ on $X$ is a principal homogeneous space for the group
of $r$-torsion points of the Jacobian of $X$. Thus there are
$r^{2g}$ such isomorphism classes.

\begin{df}
A \emph{stable, $n$-pointed, $r$-spin curve of genus $g$ and type
$\bm$} is a stable, $n$-pointed curve of genus $g$ with an
$r$-spin structure of type $\bm$. For any $\bm=(m_1,\ldots,m_n)$,
the stack of connected, stable, $n$-pointed $r$-spin curves of
genus $g$ and type $\bm$  is denoted by $\mgnrmbar$.

The disjoint union $\displaystyle\coprod_{\substack{\mathbf{m} \\
0 \leq m_i <r}} \mgnrmbar$ is denoted by $\mgnrbar$.
\end{df}

\begin{rem} \label{rem:restrict}
No information is lost by restricting $\bm$ to the range $0\le m_i
\le r-1$, since when $\bm \equiv \bm' \pmod r$ every $r$-spin
curve of type $\bm$ naturally gives an $r$-spin curve of type
$\bm'$  simply by $$ \ce_d \mapsto \ce_d \tensor \co \Big( \sum
\frac{m_i-{m'}_i}{d} p_i \Big). $$ Thus $\mgnrmbar$ is canonically
isomorphic to $\overline{\scheme{M}}^{1/r,\mathbf{m}'}_{g,n}$.
Since one of our main goals in this note is to relate the virtual
classes for different values of $\bm$, including those outside the
range $0 \le m_i \le r-1$, we will always consider
$\overline{\scheme{M}}^{1/r,\mathbf{m}'}_{g,n}$ to be equal, via
this canonical isomorphism, to the appropriate component of
$\mgnrbar$.
\end{rem}

The stack $\mgnrbar$ is a smooth Deligne-Mumford stack, finite
over $\mgnbar$, with a projective, coarse moduli space.  For $g>1$
the spaces $\mgnrmbar$ are irreducible if $\gcd(r, m_1, \ldots,
m_n)$ is odd, and they are the disjoint union of two irreducible
components if $\gcd(r, m_1, \ldots, m_n)$ is even. When the genus
$g$ is zero,  the coarse moduli space
$\overline{M}^{1/r,\bm}_{0,n}$ is either empty (if $r$ does not
divide $2 + \sum m_i$), or is canonically isomorphic to the moduli
space $\mbar_{0,n}$. Note, however, that this isomorphism does not
give an isomorphism of stacks, since the automorphisms of elements
of ${\M}^{1/r,\bm}_{0,n}$ vary differently from the way that
automorphisms of the underlying curves vary.  In any case,
$\mbar^{1/r,\bm}_{0,n}$ is always irreducible.

When the genus is one, the stack $\mbar^{1/r,\mathbf{m}}_{1,n}$ is
the disjoint union of $d$ irreducible components, where $d$  is
the number of divisors of $\gcd(r,m_1, \ldots, m_n)$.

Throughout this paper we will denote the forgetful morphism by
$p:\mgnrbar \rTo \mgnbar $, and the universal curve  by $\pi:\cgnr
\rTo \mgnrbar$.  As in the case of the moduli space of stable
curves, the universal curve possesses canonical sections
$\sigma_i:\mgnrbar \rTo \cgnr$ for $i\,=\,1,\ldots,n$.  Unlike the
case of stable curves, however, the universal curve
$\scheme{C}^{1/r, \bm}_{g,n} \rTo \mgnrmbar$ is not obtained by
considering $(n+1)$-pointed $r$-spin curves. The universal curve
$\scheme{C}^{1/r,   \bm}_{g,n}$ is birationally equivalent to
$\overline{\scheme{M}}^{1/r, (m_1, m_2, \dots, m_n, 0)}_{g,n+1}$,
but they are not isomorphic.

\subsection{Natural Cohomology Classes on $\mgnrbar$}

All of the usual classes on  $\mgnbar$ pull back to classes on
$\mgnrbar$.  We will abuse notation and use the same symbol for
these classes regardless of whether they are on $\mgnbar$ or
$\mgnrbar$.  In particular, for each $i$ with $0 \le i \le n$ we
have the classes $\psi_i := \sigma_i^*(\omega)$, and for each
positive integer $j$, we have the class $\lambda_j$, which is the
degree-$j$ term in the Chern polynomial of the $K$-theoretic
pushforward\footnote{In \cite{JKV} we used the topologists'
notation $\pi_!$ for the  pushforward that we are calling $R\pi_*$
here.}  of $\omega$: $$ c_t (R\pi_* \omega)
 =  1 + \lambda_1 t + \lambda_2 t^2 + \dots.$$

In addition to the classes pulled back from $\mgnbar$, there are
other natural cohomology classes on $\mgnrbar$.  These include
classes arising from $\ce_r$, the $r$-th root in the universal
$r$-spin structure on the universal curve over $\mgnrmbar$:
$$\psit_i(\bm) := \psit_i := \sigma_i^* (\ce_r),$$ and $\mu_j$ is
the  degree-$j$ term in the Chern character of the pushforward  of
$\ce_r$: $$ ch (R\pi_* \ce_r)  =  1 + \mu_1 t + \mu_2 t^2 +
\dots.$$

In \cite{JKV}, we established some relationships between these
classes.  Especially useful for the current work is the following.

\begin{prop}[Proposition 2.2 of \cite{JKV}]

The line bundles $\sigma^*_i\omega $ and $\sigma^*_i(\ce_r)$ on
the stack $\mgnrmbar$ are related by $$ r \sigma^*_i(\ce_r) \cong
(m_i +1) \sigma^*_i (\omega).$$ Therefore, in $\pic \mgnrmbar
\otimes \mathbb{Q}$ and in $H^{2}(\mgnrmbar,\nq)$  we have
\begin{equation}\label{eq:psit-to-psi}
\psit_i(\bm)= \frac{m_i+1}{r} \psi_i. \end{equation}
\end{prop}

It is important to note that these relations hold for all choices
of $\bm = (m_1,\dots,m_n)$, with no restriction on the range of
the $m_i$.

An immediate corollary is the following.

\begin{crl} \label{cor:vanish}
If $\bm = (m_1, \dots, m_n)$, and if $m_i = -1$ for some $i$, then
we have $$\psit_i(\bm) = 0.$$
\end{crl}

\subsection{The Virtual Class and the Vanishing Axiom}

In \cite[\S4.1]{JKV} we give axioms for a so-called \emph{virtual
class} $\cv$, where for each $\bm = (m_1,\dots,m_n)$ with $0 \le
m_i \le r-1$ the virtual class $\cv(\bm)$ is a cohomology class in
$H^{2D}(\mgnrmbar,\nq)$.  Here the dimension $D$ is given by
\begin{equation} \label{eq:deg} D =
\frac{1}{r}\bigl( (r-2)(g-1)+\sum_{i=1}^n m_i \bigr).
\end{equation}
The class $\cv$ plays the role of Gromov-Witten classes (capped 
with the virtual fundamental class of the moduli space of stable 
maps) in the theory of Gromov-Witten invariants of a smooth, 
projective variety.

Although the axioms stated in \cite[\S4.1]{JKV} are only for $\bm$
in the range $0 \le m_i \le r-1$, these axioms make sense for all
non-negative choices of $m_i$, and indeed, the
constructions of $\cv$ given in \cite{JKV} 
for $g=0$ (and all $r>1$) satisfy the axioms for all choices of
non-negative $m_i$. Even better, these constructions make sense
and satisfy the axioms in the case that one (but not more) of the
$m_i$ is equal to $-1$, and the rest are non-negative.  The axioms
are generally inconsistent in the case that two or more of the
$m_i$ are negative, although there may be some special cases
                   with two or more negative $m$'s,
where the virtual class
exists.

One of the axioms that the virtual class must satisfy is the
vanishing axiom \cite[\S4.1, Axiom 4]{JKV}:

\begin{vax} \label{ax:vax}
\emph{If $\bm=(m_1,\dots, m_n)$ has at least one $i$ such that
$m_i = r-1$, and if all $m_j$ are non-negative, then $\cv(\bm)=0.$
}
\end{vax}

This axiom seems very strange, in that it has no clear counterpart
in Gromov-Witten theory, but it is a straightforward consequence
of the Descent Axiom, which is satisfied by the constructions of
\cite{JKV}
for $g=0$.

\begin{dax}\label{ax:dax}
Let $r\geq 2$ be an integer, and $\bm = (m_1,\dots,m_n)$ be an
$n$-tuple of integers such that either $m_i \geq 0$ for all
$i=1,\ldots,n$, or there exists an integer  $1\leq j\leq n$ such
that $m_i \geq 0$ for all $i\not=j$ and $m_j = -1$.
Let $\cv(\bm)$ denote the virtual class on $\M_{g,n}^{1/r,\bm}$.
Let $r \bdelta_i$ denote an $n$-tuple of integers which is $r$ in
the $i$th position and zero in all others then
\[
\cv(\bm+ r \bdelta_i) = -\psit_i(\bm) \cv(\bm)
\]
\end{dax}

An immediate result is the following proposition.

\begin{prop}
The Descent Axiom (\ref{ax:dax}) implies the Vanishing Axiom
(\ref{ax:vax}).
\end{prop}
\begin{proof}
If $\bm$ has all $m_j$ non-negative and some $m_i = r-1$, then
$\bm' := \bm - r \bdelta_i$ meets the conditions for the Descent
Axiom to apply, so we have $$ \cv(\bm) = -\psit_i(\bm-r\bdelta_i)
\cv(\bm-r\bdelta_i).$$ But by Corollary~\ref{cor:vanish}, we have
$$\psit_i(\bm-r\bdelta_i) = 0.$$
\end{proof}

\section{The Origin of Gravitational Descendants}
\label{lps}

In this section, we show that the Descent Axiom gives a geometric
explanation of the usual definition of gravitational descendants
in an $r$-spin \cft.

Let $\bm := (m_1,\ldots,m_n)$ consist of integers $m_i$ such that
$0\leq m_i \leq r-1$ for all $i=1,\ldots,n$. Let $a_i$ be
nonnegative integers for all $i=1,\ldots,n$. Recall that the
\emph{$n$-point correlators of genus $g$} are defined \cite{JKV,W}
by the  formula
\begin{equation}
\label{eq:OldGrav} \cor{\tau_{a_1,m_1}\cdots \tau_{a_n,m_n}}_g :=
r^{1-g} \int_{[\M_{g,n}^{1/r}]} \cv(\bm) \prod_{i=1}^n
\psi_i^{a_i}.
\end{equation}
These correlators assemble into \emph{the large phase space
potential function} $\Phi(\bt)$ in $\lambda^{-2}
\nq[[\bt,\lambda^2]]$, where $\Phi(\bt) := \sum_{g=0}^\infty
\Phi_g(\bt) \lambda^{2g-2}$ and
\[
\Phi_g(\bt) := \cor{\exp(\bt\cdot\btau)}_g,
\]
where $\lambda$ and $\bt := \{ t_a^m \}$ (for nonnegative integers
$a$ and integers $0\leq m\leq r-1$) are formal parameters.
Furthermore, we define $\bt\cdot\btau :=
\sum_{a=0}^\infty\sum_{m=0}^{r-1} t_a^m \tau_{a,m}.$

The \emph{small phase space potential} is the generating function
$\chi(\bx) :=\Phi(\bt)$ where $\bx := (x^1,\ldots,x^{r-1})$ and,
on the right hand side, we have set $x^m := t_0^m $ and $t_a^m :=
0$ for all $a > 0$.

The function $\chi(\bx)$ is a generating function for correlators
which do not contain any $\psi$ classes. Such correlators are the
analogs, in the $r$-spin \cft\  of the Gromov-Witten invariants of
a smooth, projective variety. Correlators where $a_i$ is nonzero
for some $i$ are called \emph{the gravitational descendants}.

Unlike the case of Gromov-Witten invariants, however, the only
nonzero terms in $\chi(\bx)$ come from intersection numbers on
$\M_{0,n}^{1/r}$, for dimensional reasons (although higher
contributions are present in $\Phi(\bt)$). In some ways, the
potential $\chi(\bx)$ behaves as though the $r$-spin \cft\
corresponded to Gromov-Witten invariants of a variety with
fractional dimensional cohomology classes \cite{JKV,JKV2}.

As mentioned in the introduction, the appearance of the $\psi$
classes in the definition of the gravitational descendants seems
rather mysterious from a purely geometric perspective. It is also
unnatural to restrict only to intersection numbers of $\cv(\bm)$
on $\M_{g,n}^{1/r,\bm}$ corresponding to $n$-tuples of nonnegative
integers $\bm = (m_1,\ldots,m_n)$ where $m_i\leq  r-1$ for all
$i=1,\ldots,n$.

On the other hand, it is geometrically more natural to consider
the correlators
\begin{equation}
\label{eq:NewGrav} \cort{\taut_{\mt_1}\cdots\taut_{\mt_n}}_g :=
r^{1-g} \int_{[\M_{g,n}^{1/r}]} \cv(\bmt)
\end{equation}
where $\bmt := (\mt_1,\ldots,\mt_n)$ are any $n$-tuple of
nonnegative integers. They assemble into a potential function
$\chit(\bTt) := \sum_{g=0}^\infty \chit_g(\bTt) \lambda^{2g-2}$ in
$\lambda^{-2} \nq[[\bTt,\lambda^2]]$ where
\[
\chit_g(\bTt) := \cort{\exp(\bTt\cdot\btaut)}_g.
\]
Here $\bTt := \{ \Tt^\mt \}_{\mt\geq 0}$ and $\lambda$ are formal
parameters and we  define $\bTt\cdot\btaut = \sum_{\mt\geq 0}
\Tt^\mt \taut_\mt$. If one sets $x^m := \Tt^m$ for all $0\leq
m\leq r-1$, and if one sets $\Tt^\mt :=0$ otherwise, then the
function $\chit(\bTt)$ clearly is equal to $\chi(\bx)$ where $\bx
= (x^1,\ldots, x^{r-1})$.  Remarkably, the functions $\chit(\bTt)$
and $\Phi(\bt)$ are the same.

\begin{prop}
Let $\ba := (a_1,\ldots,a_n)$ be an $n$-tuple of nonnegative
integers and $\bm := (m_1,\ldots,m_n)$ be integers such that
$0\leq m_i \leq r-1$ for all $i=1,\ldots,n$.
           Let $\bmt = (\mt_1,\dots,\mt_n) = r\ba + \bm.$
The following equation holds
\[
           \cv(\bmt) = \cv(\bm) \prod_{i=1}^n
\left(\left(\frac{-\psi_i}{r}\right)^{a_i}
[r(a_i-1)+m_i+1]_r\right)
\]
where for all $0\leq m\leq r-1$,
\[
[r(a-1)+m+1]_r := \prod_{i=1}^{a} (r(a-i)+m+1)
\]
if $a\geq 1$. If $a=0$ then we define $[r(0-1)+m+1]_r := 1$.
Furthermore, $\cv(\bmt) = 0$
if, for some $i=1,\ldots,n$, we have that $\mt_i = a_i r - 1$ for
some nonnegative $a_i$.
\end{prop}

\begin{proof}
This follows from repeated application of the Descent Axiom and
Equation \ref{eq:psit-to-psi}.
\end{proof}

\begin{crl}
Let $\bt := \{ t_a^m \}$ where $a\geq 0$ and $0\leq m\leq r-1$.
Let $\bTt := \{ \Tt^\mt \}$ where $\mt$ runs over the nonnegative
integers then
\[
\chit(\bTt) = \Phi(\bt)
\]
under the identification
\[
\Tt^{a r + m} = \frac{(-1)^a r^a}{[r(a-1)+m+1]_r} t_a^m
\]
where $a$ and $m$ are integers such that $a\geq 0$ and $0\leq
m\leq r-1$.
\end{crl}

\begin{rem}
This Corollary suggests that the proper geometric definition of
gravitational descendants should be given by Equation
\ref{eq:NewGrav} rather than the usual definition given in
Equation \ref{eq:OldGrav}.
\end{rem}

It is also amusing to note that the equations of the $\KdV_r$
($r$-th Gelfand-Dickey) hierarchy \cite{W,JKV} simplify somewhat
when written in terms of the variables $\Tt^{a r + m}$ instead of
the $t_a^m$. The $\KdV_r$ equations are, for all nonnegative
integers $a$ and $m$,
\[
i \sqrt{r} \left( a + \frac{m+1}{r} \right) \frac{\partial Q}{\partial
\Tt^{a r + m}} = [Q_+^{a + \frac{m+1}{r}},Q]
\]
where
\[
Q := D^r - \sum_{m=0}^{r-2} u_m(x) D^m
\]
and $D = \frac{i}{\sqrt{r}} \frac{\partial}{\partial x}$.

Whereas, the usual equations, written in the $t_a^m$, are
\[
i [ a r + m+1]_r  \frac{\partial Q}{\partial t_a^m} =
(-1)^a r^{a+(1/2)} [Q_+^{a + \frac{m+1}{r}},Q].
\]

\section{The Convex Case}\label{sec:zero}

\

In this section, we prove that the descent property of the
$r$-spin virtual class is satisfied in the convex case.

Recall that an $r$-spin structure on a family of curves $\pi : X
\to T$ is said to be {\em convex\/} if the $r$-th root sheaf
$\ce_r$ satisfies $$\pi_* \ce_r =0.$$ In particular, on a family
of irreducible curves of genus $g$ with $n$ marked points, an
$r$-spin structure of type $\bm$ with
$|\bm|:=\sum_{i=1}^n m_i > 2g-2$ is always convex since $\ce_r$ in
this case has negative degree on each fiber $X_t$. (This is not
true for families with reducible curves because they may have
components of                     positive
genus
           on which $|\bm|$ is too small.)       
\begin{thm}
If the universal $r$-spin structure on a (connected) component
of $\mgnrmbar$ is convex then the Descent Axiom (\ref{ax:dax})
holds.
\end{thm}

\begin{proof}
The convexity axiom of an $r$-spin virtual class (cf.~\cite{JKV})
gives that
$$ \cv(\bm) = c_D(R\pi_*\ce_r)=(-1)^D c_D(R^1\pi_*\ce_r), $$ where
$$ D = \frac{1}{r}\left((r-2)(g-1)+\sum_{i=1}^n m_i\right), $$ and
$\ce_r$ is the $r$-th root sheaf of the universal $r$-spin
structure.

The $r$-th root sheaf $\ce'_r$ of the universal $r$-spin structure
of type $\bm'=\bm+r \bdelta_i$ is isomorphic to
$\ce_r(-p_i)=\ce_r\tensor \co(-p_i)$ under the identification
$$\mgnbar^{1/r,\bm'}=\mgnrmbar.$$ Since $\ce'_r$ is a subsheaf of
$\ce_r$, it is also convex and $$ D'=\rk R^1\pi_*\ce'_r=1+\rk
R^1\pi_*\ce_r=D+1. $$ Therefore,
\begin{eqnarray*}
\cv(\bm+r\bdelta_i)&=&(-1)^{D'}c_{D'}(R^1\pi_*\ce_r(-p_i))\\
&=&-(-1)^Dc_D(R^1\pi_*\ce_r)c_1(p_i^*\ce_r)=
-\cv(\bm)\psit_i(\bm),
\end{eqnarray*}
where, in the second equality, we have used the  following simple
lemma.
\end{proof}

\begin{lm}
Let $\pi: X\to T$ be a flat family of curves, $q:T\to X$ be a
section of $\pi$, and  $\ce$ be a coherent sheaf on $X$. Assume 
that the image $q(T)$ of $q$ is disjoint from the singular locus 
of $\pi$ and the restriction of $\ce$ to $q(T)$ is locally free of 
rank one. 

If $\pi_*\ce=0$ then $$ c_{top}(R^1\pi_*\ce') = c_1(q^*\ce)
c_{top}(R^1\pi_*\ce), $$ where $\ce'=\ce(-q)$ and $c_{top}$
denotes the top Chern class of a vector bundle.
\end{lm}

\begin{proof}

Consider the exact sequence of sheaves on $X$ $$ 0\to \ce' \to \ce
\to \ce|_q \to 0~. $$ Since $\pi_*\ce=0$ and  $R^1\pi_*(\ce|_q) =
0$, the corresponding long exact sequence for the functor $R\pi_*$
gives a short exact sequence of bundles on $T$ $$ 0 \to q^*\ce \to
R^1\pi_*\ce' \to R^1\pi_*\ce \to 0. $$ Now the statement of the
lemma follows from multiplicativity of the total Chern class.
\end{proof}

In~\cite[Proposition 4.4]{JKV} we proved that in the genus zero
case, the universal $r$-spin structure is convex for all $r$ and
$\bm$, such that $m_i\ge -1$ for all $i$ and at most one $m_i$
equal to $-1$. This yields the following.
\begin{crl}
The Descent Axiom (\ref{ax:dax}) holds in the case $g=0$.
\end{crl}

\begin{rem}
In~\cite{PoVa}, a class $c$ on $\M_{g,n}^{1/r,\bm}$ was 
constructed in the cases where $\bm = (m_1,\dots,m_n)$ has all 
$m_i$ nonnegative. The class $c$ was shown to have the correct 
dimension and some of the axioms of \cite{JKV} for a spin virtual 
class were verified. Furthermore, $c$ was shown to satisfy the 
Descent Axiom whenever it is defined. However, this does not imply 
the Vanishing Axiom, since the class $c$ is not defined when any 
$m_i$ is negative.
\end{rem}

\bibliographystyle{amsplain}

\providecommand{\bysame}{\leavevmode\hbox
to3em{\hrulefill}\thinspace}

\end{document}